\documentclass[12pt]{article}
\usepackage{amssymb}
\usepackage{latexsym,bm}
\usepackage{graphicx}
\usepackage{amsmath}
\usepackage{mathrsfs}

\date{}
\begin{document}
\title{On almost self-centered graphs and almost peripheral graphs\footnote{E-mail addresses:
{\tt huyanan530@163.com}(Y.Hu),
{\tt zhan@math.ecnu.edu.cn}(X.Zhan).}}
\author{\hskip -10mm Yanan Hu and Xingzhi Zhan\thanks{Corresponding author.}\\
{\hskip -10mm \small Department of Mathematics, East China Normal University, Shanghai 200241, China}}\maketitle
\begin{abstract}
 An almost self-centered graph is a connected graph of order $n$ with exactly $n-2$ central vertices, and an almost peripheral graph is a connected graph of order $n$  with exactly $n-1$ peripheral vertices. We determine (1) the maximum girth of an almost self-centered graph of order $n;$ (2) the maximum independence number of an almost self-centered graph of order $n$ and radius $r;$ (3) the minimum order of a $k$-regular almost self-centered graph  and (4) the maximum size of an almost peripheral graph of order $n;$ (5) which numbers are possible for the maximum degree of an almost peripheral graph of order $n;$ (6) the maximum number of vertices of maximum degree in an almost peripheral graph of order $n$ whose maximum degree is the second largest possible. Whenever the extremal graphs have a neat form, we also describe them.

\end{abstract}



\section{Introduction}

We consider finite simple graphs. The {\it order} of a graph is its number of vertices, and the {\it size} its number of edges.  We denote by $V(G)$ and $E(G)$ the vertex set and edge set of a graph $G$ respectively. Denote by $d_{G}(u,v)$ the distance between two vertices $u$ and $v$ in $G.$  The {\it eccentricity}, denoted by $ecc_G(v),$ of a vertex $v$ in a graph $G$ is the distance to a vertex farthest from $v.$ Thus $ecc_G(v)={\rm max}\{d_G(v,u)| u\in V(G)\}.$
If the graph $G$ is clear from the context, we omit the subscript $G.$ If $ecc(v)=d(v,x),$ then the vertex $x$ is called an {\it eccentric vertex} of $v.$
The {\it radius} of a graph $G,$ denoted ${\rm rad}(G),$ is the minimum eccentricity of all the vertices in $V(G),$ whereas the {\it diameter}
of $G,$ denoted ${\rm diam}(G),$ is the maximum eccentricity. A vertex $v$ is a {\it central vertex} of $G$ if $ecc(v)={\rm rad}(G).$ The {\it center} of a graph $G,$ denoted $C(G),$ is the set of all central vertices of $G.$ A vertex $u$ is a {\it peripheral vertex} of $G$ if $ecc(u)={\rm diam}(G).$ 
The {\it periphery} of $G$ is the set of all peripheral vertices of $G.$ A graph with a finite radius or diameter is necessarily connected.

If ${\rm rad}(G)={\rm diam}(G),$ then the graph $G$ is called {\it self-centered.} Thus, a self-centered graph is a graph in which every vertex is a central vertex. This class of graphs have been extensively studied. See [2] and the references therein. Since a nontrivial graph has at least two peripheral vertices,
a connected non-self-centered graph of order $n$ has at most $n-2$ central vertices. The following concept was introduced in [4].

{\bf Definition 1.} A connected graph of order $n$ is called {\it almost self-centered} if it has exactly $n-2$ central vertices.

Since every graph has at least one central vertex, a connected graph of order $n$ has at most $n-1$ peripheral vertices. The following concept was introduced in [5].

{\bf Definition 2.} A connected graph of order $n$ is called {\it almost peripheral} if it has exactly $n-1$ peripheral vertices.

In this paper we investigate several extremal problems on these two classes of graphs. In particular, we determine (1) the maximum girth of an almost self-centered graph of order $n;$ (2) the maximum independence number of an almost self-centered graph of order $n$ and radius $r;$ (3) the minimum order of a $k$-regular almost self-centered graph  and (4) the maximum size of an almost peripheral graph of order $n;$ (5) which numbers are possible for the maximum degree of an almost peripheral graph of order $n;$ (6) the maximum number of vertices of maximum degree in an almost peripheral graph of order $n$ whose maximum degree is the second largest possible. Whenever the extremal graphs have a neat form, we also describe them.

In Section 2 we treat almost self-centered graphs, and in Section 3 we treat almost peripheral graphs.

For graphs $G$ and $H,$ the notation $G+H$ means the disjoint union of $G$ and $H.$ 
A {\it dominating vertex} in a graph of order $n$ is a vertex of degree $n-1.$ Two vertices $u$ and $v$ on a cycle $C$ of length $n$ are called {\it antipodal vertices} if $d_C(u,v)=\lfloor n/2\rfloor.$ An {\it $(x,y)$-path} is a path with endpoints $x$ and $y.$  A {\it diametral path} in a graph $G$ is a shortest  $(x,y)$-path of length ${\rm diam}(G).$ We list some notations which will be used:
\newline\indent $C_n$:  the cycle of order $n,$  $P_n$: the path of order $n,$ $K_n$: the complete graph of order $n,$
\newline\indent $\overline{G}$: the complement of the graph $G,$ \,\, $e(G)$: the size of the graph $G,$
\newline\indent $\delta(G)$: the minimum degree of vertices of  the graph $G,$
\newline\indent $\Delta(G)$: the maximum degree of vertices of  the graph $G,$
\newline\indent $\alpha(G)$: the independence number of the graph $G,$ \,\, $g(G)$: the girth of the graph $G,$
\newline\indent ${\rm deg}(v)$: the degree of the vertex $v,$ \,\, $N(v)$: the neighborhood of the vertex $v,$
\newline\indent $N[v]$: the closed neighborhood of the vertex $v;$ i.e., $N[v]=N(v)\cup \{v\},$
\newline\indent $N_i(v)$: the $i$-th neighborhood of the vertex $v;$ i.e.,  $N_i(v)=\{x\in V(G)|\, d(v,x)=i\}.$

It is known [6, p.288] that if $G$ is a connected graph satisfying ${\rm diam}(G)\ge {\rm rad}(G)+2,$ then every integer $k$ with
${\rm rad}(G)<k<{\rm diam}(G)$ is the eccentricity of some vertex. Thus if $G$ is an almost self-centered graph or an almost peripheral graph,
then the vertices of $G$ have only two distinct eccentricities and hence ${\rm diam}(G)={\rm rad}(G)+1.$

\section{Almost self-centered graphs}

A {\it binocle} is a graph that consists of two cycles $C,\,D$ and a $(u,v)$-path $P$ such that $V(P)\cap V(C)=\{u\},$ $V(P)\cap V(D)=\{v\}$
and $V(C)\cap V(D)\subseteq V(P).$  Here we allow the possibility that $P$ has length $0;$ i.e., $P$ is a vertex. Note also that if $P$ is nontrivial,
then $C$ and $D$ are vertex-disjoint. A {\it theta} (or {\it theta graph}) is a graph that consists of three internally
vertex-disjoint paths sharing the same two endpoints. $\theta_{a,b,c}$ will denote the theta consisting of three paths with lengths $a,$ $b$ and $c$ respectively. A binocle and a theta is depicted in Figure 1.

\vskip 3mm
\par
 \centerline{\includegraphics[width=4.5in]{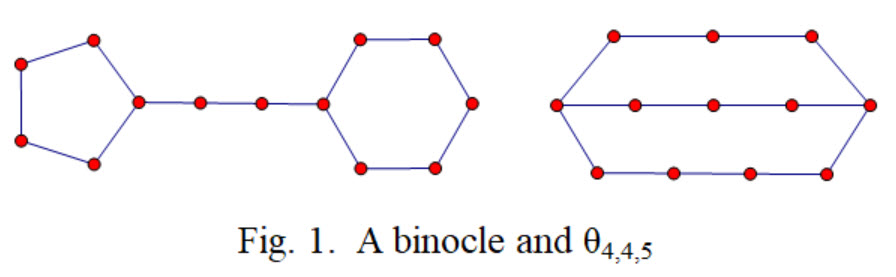}}
\par

We make the convention that the girth of an acyclic graph is undefined. Thus whenever we talk about the girth of a graph, the graph is not acyclic.
A connected graph is said to be {\it unicyclic} if it contains exactly one cycle. Recall that a connected graph of order $n$ is unicyclic if and only if
it has size $n$ [7, p.77].

{\bf Lemma 1.} {\it Let $G$ be a unicyclic graph of order $n\ge 6.$ Then $G$ is almost self-centered if and only if $n$ is odd and $G$ is the graph
obtained from $C_{n-1}$ by attaching one edge.}

{\bf Proof.} Suppose $G$ is almost self-centered. We have $e(G)=n$ and $G\neq C_n.$ It is known [3] that the center of any connected
graph lies within one block. Let $B$ be the block of $G$ in which $C(G)$ lies. Then $B$ is unicyclic and $\delta(B)=2.$ Thus $B$ is a cycle of order
$n-2$ or $n-1.$ If $B=C_{n-2},$ let $V(G)\setminus V(B)=\{x,y\}.$ Then $x$ and $y$ are leaves. Since $x$ and $y$ are the only two peripheral vertices,
their neighbors are a pair of antipodal vertices of the cycle $B.$ But then $B$ contains a vertex whose eccentricity in $G$ is ${\rm diam}(G)-2,$
contradicting the assumption that $G$ is almost self-centered. Hence $B=C_{n-1}$ and $G$ is the graph obtained from $C_{n-1}$ by attaching one edge.
Since $G$ is almost self-centered, $n$ is odd.

Conversely, it is easy to verify that this graph is almost self-centered. \hfill $\Box$

{\bf Lemma 2.} {\it If $G$ is a connected graph of order $n$ and size $n+1$ with $\delta(G)=2,$ then $G$ is either a binocle or a theta.}

{\bf Proof.} Since $\sum_{x\in V(G)}{\rm deg}(x)=2n+2$ and $\delta(G)=2,$ the degree sequence of $G$ is $(2,2,\ldots,2,4)$ or $(2,2,\ldots,2,3,3).$
In the former case, $G$ is a graph consisting of two cycles sharing a common vertex, which is a binocle, while in the latter case, $G$ is a theta.
\hfill $\Box$

Lemma 2 can also be proved easily using induction on the order.

{\bf Lemma 3.} {\it Let $a,b,c$ be positive integers with $a\le b\le c.$ Then ${\rm rad}(\theta_{a,b,c})=\lfloor (a+c)/2\rfloor$ and
${\rm diam}(\theta_{a,b,c})=\lfloor (b+c)/2\rfloor.$ Consequently $\theta_{a,b,c}$ is self-centered if and only if $b=a$ if $a+c$ is odd
and $b\le a+1$ if $a+c$ is even. Also ${\rm diam}(\theta_{a,b,c})={\rm rad}(\theta_{a,b,c})+1$ if and only if $a+1\le b\le a+2$ if $a+c$
is odd and $a+2\le b\le a+3$ if $a+c$ is even. }

{\bf Proof.} Easy verification. \hfill $\Box$

{\bf Lemma 4.} {\it Let $G$ be a connected graph of order $n$ and size $n+1$ with $\delta (G)=2.$ Then $G$ is almost self-centered if and only
if $n$ is even and $G=\theta_{1,2,n-2}.$ }

{\bf Proof.} By Lemma 2, $G$ is either a binocle or a theta.
Suppose that $G$ is almost self-centered. It is easy to see that an almost self-centered graph with minimum degree $2$ is $2$-connected.
Since a binocle has connectivity $1,$ we deduce that $G$ is a theta. Let $G=\theta_{a,b,c}$ with $a\le b\le c.$ Since $G$ is almost self-centered,
${\rm diam}(G)={\rm rad}(G)+1.$ By Lemma 3, $a+1\le b\le a+2$ if $a+c$ is odd and $a+2\le b\le a+3$ if $a+c$ is even.
 First suppose that $a+c$ is odd. We assert that $a=1.$ To the contrary, assume $a\ge 2.$ Then $b\ge a+1\ge 3.$ Let $G=\theta_{a,b,c}$  consist
 of the three $(x,y)$-paths $P_1,\, P_2,\, P_3$ of lengths $a,\, b,\, c$ respectively. Denote $r={\rm rad}(G)$ and $d=r+1={\rm diam}(G).$
 Note that $x$ and $y$ are central vertices of $G;$ i.e., $ecc(x)=ecc(y)=r.$ Let $w_1$ be the neighbor of $x$ on $P_2$ and let $w_2$ be the neighbor
 of $y$ on $P_2.$ Let $x_1$ and $x_2$ be the two antipodal vertices of $x$ on the odd cycle $C=P_1\cup P_3$ where $d_C(x_2,y)=d_C(x_1,y)+1.$
 Then $d_G(w_1,x_2)\ge r+1=d.$ Thus both $w_1$ and $x_2$ are peripheral vertices. Similarly, $w_2$ is a peripheral vertex. But then $G$ contains
 at least three peripheral vertices, a contradiction.

 The case when $a+c$ is even can be treated similarly. Hence $a=1.$ Lemma 3 implies that $b\ge 2$ if $a+c$ is odd and $b\ge 3$ if $a+c$ is even.
 If $b\ge 3,$ using the above argument we obtain contradictions. Thus $a+c=1+c$ is odd and $b=2.$ It follows that $n=a+b+c-1=1+(1+c)$ is even
 and $G=\theta_{1,2,n-2}.$

 Conversely, it is easy to verify that if $n$ is even then the theta $\theta_{1,2,n-2}$ is almost self-centered. \hfill $\Box$

 Now we are ready to state and prove the first main result.

{\bf Theorem 5.} {\it Let $g(n)$ denote the maximum girth of an almost self-centered graph of order $n$ with $n\ge 5.$ Then
$$
g(n)=\begin{cases} n-1 \quad {\rm if}\,\,\,n\,\,\,{\rm is}\,\,\,{\rm odd},\\
4\lfloor n/6\rfloor \quad {\rm if}\,\,\,n\,\,\,{\rm is}\,\,\,{\rm even}\,\,\,{\rm and}\,\,\,n\neq 10,\\
5 \quad {\rm if}\,\,\,n=10.
\end{cases}
$$
Furthermore, if $n\ge 12$ and $6$ divides $n,$ then $g(n)$ is attained uniquely by the graph obtained from $\theta_{n/3,n/3,n/3}$ by attaching an edge
to a vertex of degree three. }

{\bf Proof.} Let $G$ be an almost self-centered graph of order $n\ge 5.$ Clearly $G\neq C_n.$ Hence $g(n)\le n-1.$ On the other hand, if $n$ is odd,
then the graph obtained from $C_{n-1}$ by attaching an edge is almost self-centered and has girth $n-1.$ Hence $g(n)=n-1$ if $n$ is odd.

Now suppose that $n$ is even. Note that adding edges to a graph does not increase its girth. The cases $n\le 16$ can be verified by a computer search.
Using Lemma 1 and the fact [1, p.195] that a graph of order $n$ and size $n+3$ has girth at most $\lfloor 4(n+3)/9\rfloor,$ we need only
check the sizes $n+1$ and $n+2$ for a graph of order $n\le 16.$

Next suppose that $n$ is even and $n\ge 18.$ We first show that  $g(G)\le 4\lfloor n/6\rfloor.$ It is known [1, p.195] that a graph of
order $n$ and size $n+2$ has girth at most $\lfloor n/2\rfloor +1.$ The inequality $\lfloor n/2\rfloor +1\le 4\lfloor n/6\rfloor$ for $n\ge 18$
implies that if $e(G)\ge n+2,$ then $g(G)\le 4\lfloor n/6\rfloor.$ Also, Lemma 1 excludes the possibility that $e(G)=n.$ It remains to consider the case
when $e(G)=n+1,$ and from now on we make this assumption. It is known [3] that the center of any connected graph lies within one block. Let $B$ be the block
of $G$ in which $C(G)$ lies. Since $|C(G)|=n-2$ and $e(G)=n+1,$ the size of $B$ equals its order plus one. Since $n\ge 18,$ $B$ is $2$-connected and
$\delta(B)=2.$ By Lemma 2, $B$ is a theta. Let $B=\theta_{a,b,c},$ which consists of three $(x,y)$-paths $P_1,\,P_2,\,P_3$ whose lengths are $a,\,b,\,c$
respectively with $a\le b\le c.$

Since the eccentricities of two adjacent vertices differ by at most one, every leaf of $G$ is a peripheral vertex. Hence $G$ has at most two leaves.
We first exclude the possibility of two leaves. To the contrary, assume that $G$ has two distinct leaves $u$ and $v$ whose neighbors are $s$ and $t$
respectively.

Denote $d={\rm diam}(B)$ and $f={\rm diam}(G).$ Then $d+1\le f\le d+2.$ Clearly it is impossible that $f<d.$ It is also impossible that $f=d,$
since otherwise $G$ would have at least four peripheral vertices, a contradiction. Hence $f\ge d+1.$ The inequality $f\le d+2$ follows from the fact
that adding two leaves to $B$ can increase its diameter by at most $2.$ We distinguish two cases.

Case 1. $f=d+2.$ In this case ${\rm rad}(G)=d+1.$ Since adding leaves to $B$ can increase the eccentricity of any vertex of $B$ by at most $1,$
we deduce that $B$ is self-centered. Clearly $d_B(s,t)=d\ge 2.$ Let $w$ be an internal vertex on a shortest $(s,t)$-path in $B.$ Then $ecc_G(w)=d,$
contradicting ${\rm rad}(G)=d+1.$

Case 2. $f=d+1.$ Since ${\rm rad}(G)=d,$ we have ${\rm rad}(B)\ge d-1.$ We further consider two subcases.

Subcase 2.1. ${\rm rad}(B)=d;$ i.e., $B$ is self-centered. Let $s^{\prime}$ be an eccentric vertex of $s$ in $B.$ Then $d_B(s^{\prime},s)=d,$ implying
that $d_G(s^{\prime},u)=d+1=f.$ But then $G$ has at least three peripheral vertices $u,\,v,\, s^{\prime},$ a contradiction.

Subcase 2.2. ${\rm rad}(B)=d-1.$ By Lemma 3, $a+1\le b\le a+2$ if $a+c$ is odd and $a+2\le b\le a+3$ if $a+c$ is even. It suffices to consider the two
cases: $b\ge a+2;$ $b=a+1$ and $a+c$ is odd.

First suppose $b\ge a+2.$ Since $ecc_B(x)={\rm rad}(B)=d-1$ and ${\rm rad}(G)=d,$ one of the two leaves, say $u,$ must be an eccentric vertex of $x$
in $G.$ Let $p$ be the neighbor of $x$ on $P_2.$ Using the structure of the theta $B$ and the condition $b\ge a+2,$ we deduce that $d_G(p,u)\ge d+1.$
Consequently $G$ has at least three peripheral vertices $u,\,v,\, p,$ which is a contradiction.

Next suppose that $b=a+1$ and $a+c$ is odd. If in $G,$ $x$ and $y$ have a common eccentric vertex (which must be one of the two leaves), say $u,$ then
$s$ is a common eccentric vertex of $x$ and $y$ in $B.$ Note that now $s$ lies in $P_3.$ Such a situation occurs only if $a=1,$ and hence $b=2.$
Let $q$ be the internal vertex of $P_2.$ Then $d_G(q,u)=d+1.$ Thus $G$ has at least three peripheral vertices $u,\,v,\, q,$ which is a contradiction.

If in $G,$ $x$ and $y$ do not have a common eccentric vertex, then one of $u$ and $v$ is an eccentric vertex of $x$ and the other is an eccentric vertex
of $y.$ The conditions $a+b+c=n-1,$ $b=a+1,$ and $n\ge 18$ imply that $a+c\ge 8.$ Hence $d-1=\lfloor (a+c)/2\rfloor\ge 4.$ Since $d_G(u,v)=d+1,$ we have
$d_B(s,t)=d-1\ge 4.$ Note that $s$ and $t$ lie in $P_3.$ Choose two adjacent vertices $v_1$ and $v_2$ on $P_3$ between $s$ and $t.$ Since the cycle
$P_1\cup P_3$ is odd, $v_1$ and $v_2$ have a common antipodal vertex $z$ on $P_1.$ It is easy to verify that $ecc_G(z)=d-1,$ contradicting the fact that
${\rm rad}(G)=d.$

If $G$ has no leaf, by Lemma 4 $G=\theta_{1,2,n-2}.$ Thus $g(G)=3< 4\lfloor n/6\rfloor.$

Finally we consider the case when $G$ has exactly one leaf. Let $u$ be the leaf and let $s$ be its neighbor. Note that $s\in V(B),$ since otherwise
the vertices of $G$ would have at least three distinct eccentricities, contradicting the assumption that $G$ is almost self-centered. We continue using
the notations $d={\rm diam}(B)$ and $f={\rm diam}(G).$ If $f=d,$ then $G$ would have at least three peripheral vertices, a contradiction. It is also
impossible that $f\ge d+2,$ since adding a leaf increases the eccentricity of any vertex by at most $1.$ Hence $f=d+1.$ Clearly ${\rm rad}(B)\ge d-1.$

We assert that $B$ is self-centered; i.e., ${\rm rad}(B)=d.$ To the contrary, suppose ${\rm rad}(B)=d-1.$ Then $ecc_B(x)=ecc_B(y)=d-1.$ Since
${\rm rad}(G)=d,$ we deduce that $u$ must be the common eccentric vertex of $x$ and $y,$ implying that $s$ is a common antipodal vertex of $x$ and $y$
on the cycle $P_1\cup P_3.$ As argued above, $a+c$ is odd and $a=1.$ By Lemma 3, $b=2$ or $b=3.$ If $b=2,$ let $z$ be a neighbor of $s$ on $P_3.$
Then $ecc_G(z)=d-1,$ contradicting the fact that ${\rm rad}(G)=d.$ If $b=3,$ let $w_1$ be the neighbor of $x$ on $P_2$ and let $w_2$ be the neighbor of
$y$ on $P_2.$ Then it is easy to check that $G$ has at least three peripheral vertices $u,\,w_1,\, w_2,$ which is a contradiction again. Thus $B$ is self-centered.

By Lemma 3, $b=a$ if $a+c$ is odd and $b\le a+1$ if $a+c$ is even. We have $a+b+c=n,$ and clearly $g(G)=a+b.$ There are two possibilities:
(1) $b=a$; (2) $b=a+1$ and $a+c$ is even. Denote $k=\lfloor n/6\rfloor.$

Suppose $b=a.$ We have $3a\le n,$ implying that $a\le n/3.$ If $n=6k$ or $n=6k+2,$ we obtain $g(G)=2a\le 4k.$ If $n=6k+4,$ we have $a\le 2k+1.$
The case $a=2k+1$ will be excluded. Assume $a=2k+1.$ Then $b=2k+1$ and $c=2k+2.$ Thus $a+c=b+c$ is odd. But then $G$ has at least three peripheral vertices,
a contradiction. Hence $a\le 2k$ and $g(G)=2a\le 4k.$

Suppose $b=a+1$ and $a+c$ is even. We have $a\le (n-1)/3\le 2k+1,$ where in the second inequality we have used $n\le 6k+4.$ But it is impossible that
$a=2k+1,$ since otherwise $c=2k+1<2k+2=b,$ contradicting our assumption that $b\le c.$ The conditions $a+b+c=n,$ $b=a+1$ and that both $n$ and $a+c$
are even imply that $a$ is odd. Thus $a=2k$ is also impossible. It follows that $a\le 2k-1$ and consequently $g(G)=a+b=2a+1\le 4k-1.$

Finally we prove that the upper bound $4\lfloor n/6\rfloor$ can be attained and when $6$ divides $n,$ the extremal graph is unique. Denote
$k=\lfloor n/6\rfloor.$ Let $G$ be the graph obtained from $\theta_{2k,2k,n-4k}$ by attaching an edge to one of the two vertices of degree three.
Then $G$ is an almost self-centered graph of order $n$ with girth $4\lfloor n/6\rfloor.$

Suppose $G$ is an almost self-centered graph of order $n=6k\ge 18$ with girth $4k.$ Then the above analysis shows that $G$ is a graph obtained from
$\theta_{a,b,c}$ by attaching an edge where $a=b.$ Since $g(G)=a+b=2a=4k,$ we have $a=b=2k.$ The condition $a+b+c=n$ further implies $c=2k.$ Thus the theta is
$\theta_{2k,2k,2k}.$ There is only one way to attach an edge to this theta so that the resulting graph is almost self-centered; i.e., attach the edge
to a vertex of degree three. This shows that the extremal graph is unique. The proof is complete. \hfill $\Box$

One conclusion in Theorem 5 states that if $n\ge 12$ and $6$ divides $n,$ then the extremal graph for $g(n)$ is unique. We remark that if $n$ is even
with $n\ge 14$ and $6$ does not divide $n,$ then there are at least three extremal graphs for $g(n).$ This can be seen as follows. Using the notations in the proof of Theorem 5, we may attach an edge to any vertex on $P_1$ of the theta $\theta_{2k,2k,n-2k}$ to obtain an extremal graph.

Next we consider the independence number. There is only one almost self-centered graph of order $n$ and radius $1;$ i.e., the graph obtained from $K_n$
by deleting an edge.

{\bf Theorem 6.} {\it The maximum independence number of an almost self-centered graph of order $n$ and radius $r$ with $r\ge 2$ is $n-r.$ }

{\bf Proof.} Let $G$ be an almost self-centered graph of order $n$ and radius $r.$ First recall that ${\rm diam}(G)={\rm rad}(G)+1=r+1.$
Let $P$ be a diametral path of $G.$ If $r=2,$ $P$ has order $4.$ Any independent set can contain at most $2$ of the four vertices on $P.$
Thus $\alpha (G)\le n-2.$

Suppose $r\ge 3.$ Let $x$ be a central vertex of the path $P.$ Now $P$ has order at least $5$ and any vertex on $P$ is not an eccentric vertex
of $x.$  Let $y$ be an eccentric vertex of $x.$ It is known [3] that the center of any connected graph lies within one block. Let $B$ be the block of
$G$ in which $C(G)$ lies. Then $x,\,y\in V(B).$ By Menger's theorem [7, p.167], there are two internally disjoint $(x,y)$-paths $Q_1$ and $Q_2.$
Denote by $k$ the length of the cycle $D=Q_1\cup Q_2.$ Then $k\ge 2r.$ Any independent set can contain at most $\lfloor k/2\rfloor$ vertices on $D.$
Thus $\alpha (G)\le \lfloor k/2\rfloor+(n-k)=n-\lceil k/2\rceil\le n-(k/2)\le n-r.$

Conversely we construct a graph to show that the upper bound $n-r$ can be attained. Attaching an edge to the cycle $v_1,\, v_2,\ldots,v_{2r}$ at the vertex
$v_1$ we obtain a graph $H.$ Adding $n-2r-1$ new vertices to $H$ such that each of them has $v_1$ and $v_3$ as neighbors, we obtain the graph $Z(n,r).$
It is easy to see that $Z(n,r)$ is an almost self-centered graph of order $n$ and radius $r$ with  independence number $n-r.$ The graph $Z(12,4)$ is depicted
in Figure 2. \hfill $\Box$
\vskip 3mm
\par
 \centerline{\includegraphics[width=2.4in]{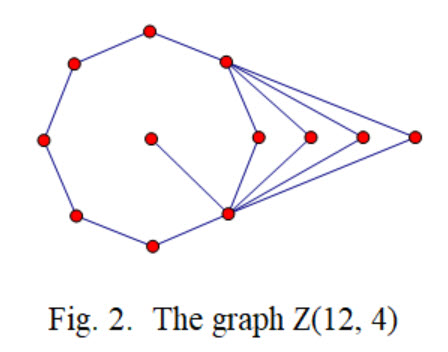}}
\par

{\bf Corollary 7.} {\it The maximum independence number of an almost self-centered graph of order $n$ with $n\ge 5$ is $n-2,$ and there are exactly two extremal
graphs.}

{\bf Proof.} By Theorem 6 and the fact that the almost self-centered graph of order $n$ and radius $1$ has independence number $2,$ we deduce that
the maximum independence number is $n-2.$

Suppose $G$ is an almost self-centered graph of order $n$ whose independence number is $n-2.$ By Theorem 6, ${\rm rad}(G)=2$ and consequently
${\rm diam}(G)=3.$ Let $P:\,x_1,x_2,x_3,x_4$ be a diametral path of $G.$ Then $x_1$ and $x_4$ are the two peripheral vertices of $G.$
Denote $S=V(G)\setminus\{x_1,x_2,x_3,x_4\}.$ $G$ has only one maximum independent set; i.e, $S\cup T$  where $T$ consists of two vertices from $P.$
There are three possible choices for $T:$ $\{x_1,\,x_3\},$ $\{x_2,\,x_4\}$ and $\{x_1,\,x_4\},$ the first two of which will yield isomorphic graphs.
Since every leaf of an almost self-centered graph is a peripheral vertex, every vertex in $S$ has degree at least $2.$ If $T=\{x_1,\,x_3\},$ then
every vertex in $S$ has $x_2$ and $x_4$ as neighbors; if $T=\{x_1,\,x_4\},$ then every vertex in $S$ has $x_2$ and $x_3$ as neighbors. Conversely,
it is easy to see that these two graphs satisfy all the requirements. \hfill $\Box$

Now we consider regular almost self-centered graphs.

{\bf Theorem 8.} {\it Let $r(k)$ denote the minimum order of a $k$-regular  almost self-centered graph. Then
$$
r(k)=\begin{cases} 12 \quad {\rm if}\,\,\,k=3,\\
2k+2 \quad {\rm if}\,\,\,k\ge 4.
\end{cases}
$$}

{\bf Proof.} Let $G$ be a $k$-regular  almost self-centered graph of order $n$, and let $x$ and $y$ be the two peripheral vertices of $G.$
There is only one almost self-centered graph of order $n$ and diameter at most $2;$ i.e., the graph obtained from $K_n$ by deleting an edge. Thus
${\rm diam}(G)\ge 3,$ implying that $N[x]\cap N[y]=\phi.$ It follows that
$$
n\ge |N[x]|+|N[y]|=(k+1)+(k+1)=2k+2.
$$

We first show $r(3)=12.$ Suppose $k=3.$ Then $n$ is even and $n\ge 2\times 3+2=8.$ We will exclude the two orders $8$ and $10.$
If $n=8,$ then ${\rm diam}(G)=3$ and $G-\{x,\,y\}$ is a $2$-regular graph of order $6,$ which must be $C_6$ or $2C_3.$ In each case, $G$ has at least four peripheral vertices, a contradiction.
  
 If $n=10,$ we deduce that ${\rm diam}(G)=3,$ since otherwise either $G$ has a vertex of degree at least $4$ or
$G$ has three peripheral vertices. Recall that $N_i(x)=\{v\in V(G)|\, d(x,\,v)=i\}.$ We have $|N_1(x)|=3$ and $|N_3(x)|=1,$ implying 
$|N_2(x)|=5.$ Here we have used the fact that $G$ has exactly two peripheral vertices. Note that each vertex in $N_2(x)$ has at least one neighbor in $N_1(x).$ Analyzing possible adjacency relations in $G-\{x,\,y\},$ we deduce that $G$ has at least four peripheral vertices, a contradiction.

Thus we have proved that $n\ge 12.$ On the other hand, the graph depicted in Figure 3 is a $3$-regular  almost self-centered graph of order $12.$
This shows $r(3)=12.$
\vskip 3mm
\par
 \centerline{\includegraphics[width=3.2in]{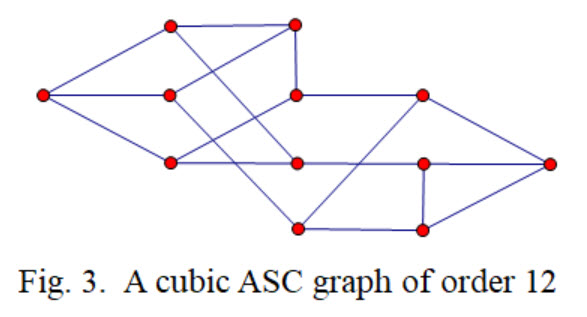}}
\par

Next suppose $k\ge 4.$ We have proved above that any $k$-regular  almost self-centered graph has  order at least $2k+2.$ To show $r(k)=2k+2,$ it suffices
to construct such a graph $R$ of order $2k+2.$ Let $V(R)=\{x_0,\,y_0\}\cup A\cup B$ where $A=\{x_1,x_2,\dots,x_k\}$ and $B=\{y_1,y_2,\dots,y_k\}.$
We use the notation $u\leftrightarrow v$ to mean that the two vertices $u$ and $v$ are adjacent. If $k$ is even, the adjacency of $R$ is defined as follows:
\begin{align*}
N(x_0)=A,\,\,\, N(y_0)=B,\,\,\, N(x_i)=\{x_{i+\frac{k}{2}},y_i,y_{i+1},\ldots,y_{i+k-3}\}\,\,\, {\rm if}\,\,\,1\le i\le k/2,\quad\quad\\
N(x_j)=\{x_{j-\frac{k}{2}},y_j,y_{j+1},\ldots,y_{j+k-3}\}\,\,\, {\rm if}\,\,\,k/2+1\le j\le k, \quad
y_i\leftrightarrow y_{i+\frac{k}{2}}, \,\,\,{\rm if}\,\,\,1\le i\le k/2.
\end{align*}
If $k$ is odd, the adjacency of $R$ is defined as follows:
\begin{align*}
N(x_0)=A,\,\,\, N(y_0)=B,\,\,\, N(x_1)=\{x_2,x_3,\ldots,x_{\frac{k+1}{2}}\}\cup \{y_1,y_2,\ldots,y_{\frac{k-1}{2}}\},\\
N(x_i)=\{x_1\}\cup (B\setminus\{y_{i-1},y_k\})\,\,\, {\rm if}\,\,\,2\le i\le (k+1)/2,\quad\quad\quad\quad\\
N(x_j)=\{x_1\}\cup (B\setminus\{y_{j-\frac{k+1}{2}}\})\,\,\, {\rm if}\,\,\,(k+3)/2\le j\le k,\quad\quad\quad\quad\\
N(y_k)=\{x_{\frac{k+3}{2}},\ldots,x_k\}\cup\{y_1,y_2,\ldots,y_{\frac{k-1}{2}}\}.\quad\quad\quad\quad\quad\quad
\end{align*}
Here the subscripts of the vertices are taken modulo $k.$ It is easy to verify that $R$ is a $k$-regular almost self-centered graph of order $2k+2$
with periphery $\{x_0,\, y_0\}$ and center $A\cup B.$  \hfill $\Box$

\section{Almost peripheral graphs}

{\bf Theorem 9.} {\it The maximum size of an almost peripheral graph of order $n$ is $\lfloor (n-1)^2/2\rfloor.$ If $n$ is odd, this maximum size
is attained uniquely by the graph $\overline{K_1+((n-1)/2)K_2};$  if $n$ is even, this maximum size is attained uniquely by the graph
$\overline{K_1+((n-4)/2)K_2+P_3}.$ }

{\bf Proof.} Use the fact that an almost peripheral graph can have at most one dominating vertex and the degree sum formula. \hfill $\Box$

In the following result we determine which numbers are possible for the maximum degree of an almost peripheral graph with a given order.

{\bf Theorem 10.} {\it There exists an almost peripheral graph of order $n\ge 7$ with maximum degree $\Delta$ if and only if
$\Delta\in \{3,\,4,\ldots,n-4,\,n-1\}.$ }

{\bf Proof.} Suppose that $G$ is an almost peripheral graph of order $n\ge 7$ with maximum degree $\Delta.$ Clearly $3\le \Delta\le n-1.$
We first exclude the two values $n-2$ and $n-3$ for $\Delta.$  To the contrary suppose $\Delta=n-2$ or $n-3.$
Note that ${\rm rad}(G)\ge 2$ and ${\rm diam}(G)={\rm rad}(G)+1\ge 3.$ Let $x\in V(G)$ with ${\rm deg}(x)=\Delta.$ There exists a vertex $y$
with $y\not\in N[x]$ such that $y$ and $x$ have a common neighbor $w.$

If $\Delta=n-2,$ then both $x$ and $w$ have eccentricity at most $2,$ implying that they are central vertices, a contradiction.

Suppose $\Delta=n-3.$ Let $z$ be the vertex outside $N[x]\cup\{y\}.$ We always have ${\rm rad}(G)=2$ and hence ${\rm diam}(G)={\rm rad}(G)+1=3.$

If $z$ and $y$ are adjacent, then $w$ is the central vertex. Since $ecc(x)=3$ and $z$ is the only possible eccentric vertex of $x,$
we deduce that $z$ is nonadjacent to any vertex in $N[x].$ It follows that $z$ is a leaf. Since $ecc(z)=3,$ $y$ is another central vertex, a contradiction.

If $z$ and $y$ are nonadjacent, then $N(x)\cap N(z)\neq\phi.$ In this case, $x$ is the central vertex. Since ${\rm diam}(G)=3,$ there
exists a $(y,z)$-path $P$ of length $2$ or $3.$ Then any internal vertex of $P$ is a central vertex different from $x,$ a contradiction.

Conversely we will show that every number in $\{3,\,4,\ldots,n-4,\,n-1\}$ can be attained. The star of order $n$ is an almost peripheral graph
with maximum degree $n-1.$ Next, for each $\Delta$ with $3\le\Delta\le n-4$ we construct an almost peripheral graph $G(n,\Delta)$ of order $n$
with maximum degree $\Delta.$ We will first construct all $G(n,3)$ for $n=7,8,\ldots,$ and then inductively construct the remaining $G(n,\Delta)$
with $\Delta\ge 4.$

$G(7,3),$ $G(8,3),$ $G(9,3),$ and $G(10,3)$ are depicted in Figure 4.
\vskip 3mm
\par
 \centerline{\includegraphics[width=4.2in]{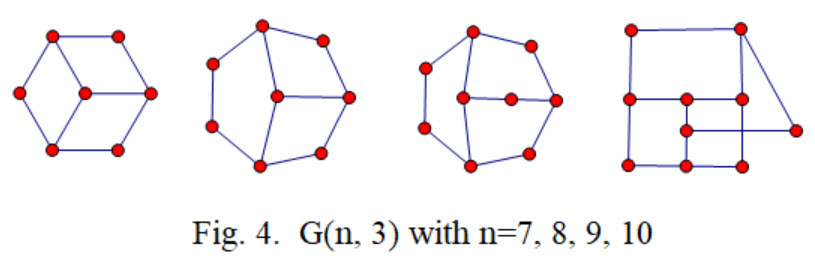}}
\par
We will need the four preliminary graphs in Figure 5.
\vskip 3mm
\par
 \centerline{\includegraphics[width=4.5in]{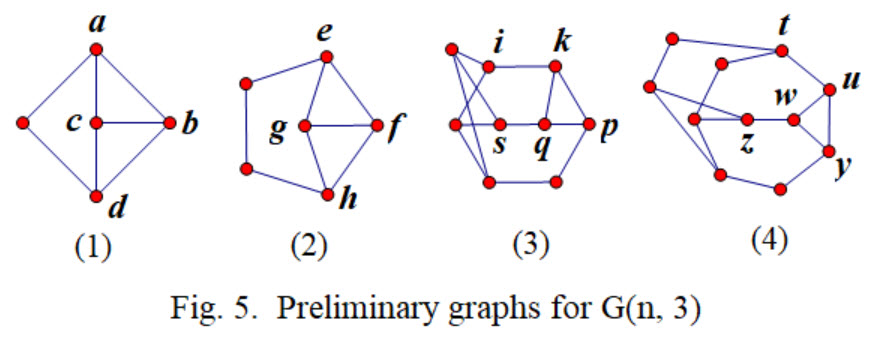}}
\par
Now let $n\ge 11$ and denote $k=\lfloor (n+5)/4\rfloor.$ If $n\equiv 3\,\,\,{\rm mod}\,\,\,4,$ $G(n,3)$ is obtained from the graph in Figure 5 (1)
by replacing the edges $ab$ and $bd$ by a path of length $k-1,$ and replacing the edges $ac$ and $cd$ by a path of length $k-2;$ if
$n\equiv 0\,\,\,{\rm mod}\,\,\,4,$ $G(n,3)$ is obtained from the graph in Figure 5 (2) by replacing the edges $ef$ and $fh$ by a path of length $k-1,$
and replacing the edges $eg$ and $gh$ by a path of length $k-2;$ if $n\equiv 1\,\,\,{\rm mod}\,\,\,4,$ $G(n,3)$ is obtained from the graph in Figure 5 (3)
by replacing the edges $ik,$ $kp,$ $sq$ and $qp$ by a path of length $k-2,$ $k-1,$ $k-3$ and $k-2$ respectively; if $n\equiv 2\,\,\,{\rm mod}\,\,\,4,$
$G(n,3)$ is obtained from the graph in Figure 5 (4) by replacing the edges $tu,$ $uy,$ $zw$ and $wy$ by a path of length $k-3,$ $k-1,$ $k-3$ and $k-2$ respectively.

For a vertex $v$ in a graph, the operation {\it duplicating $v$} means that adding a new vertex $x$ and adding edges incident to $x$ such that
$N(x)=N(v).$

Note that for $n=7,$ $n-4=3$ and that every $G(n,3)$ constructed above contains a vertex of degree $3$ that has a non-central neighbor of degree $2.$
Now suppose that we have constructed $G(n,3),$ $G(n,4),$ $\ldots,$ $G(n,n-4)$ where $G(n,\Delta)$ contains a vertex of degree $\Delta$
that has a non-central neighbor $x_{\Delta}$ of degree $2,$  $\Delta=3,\ldots, n-4,$ Then in $G(n,\Delta),$ duplicate the vertex $x_{\Delta}$ to obtain a new graph which we denote by $G(n+1,\Delta+1).$ Thus we can construct $G(n+1,3),$ $G(n+1,4),$ $\ldots,$   $G(n+1,n-3)$ which satisfy all the requirements
and the additional condition of containing a vertex of maximum degree that has a non-central neighbor of degree $2.$ Thus the inductive steps can continue.
\hfill $\Box$

Finally we consider the maximum number of vertices of maximum degree in an almost peripheral graph. {\it Blowing up a vertex $v$ in a graph into the complete graph $K_t$} is the operation of replacing $v$ by $K_t$ and adding edges joining each vertex in $N(v)$ to each vertex in $K_t.$

{\bf Definition 3.} A vertex $v$ in a graph $G$ is called a {\it top vertex} if ${\rm deg}(v)=\Delta(G).$

{\bf Theorem 11.} {\it The maximum number of top vertices in an almost peripheral graph of order $n\ge 8$ with maximum degree $n-4$
is $n-5$ and this maximum number is uniquely attained by the graph obtained from the graph of order $7$ in Figure 4 by blowing up a non-central vertex
of degree $3$ into $K_{n-6}.$ }

{\bf Proof.}  First, it is easy to verify that the extremal graph given in Theorem 11 is an almost peripheral
graph of order $n$ with maximum degree $n-4$ that has $n-5$ top vertices. Let $G$ be an almost peripheral graph of order $n\ge 8$ with maximum degree $n-4.$
We may suppose that $G$ has at least three top vertices, since otherwise the number of top vertices in $G$ is less than $n-5.$ Recall that 
${\rm diam}(G)={\rm rad}(G)+1\ge 2.$

Let $x$ be a peripheral vertex of degree $n-4.$ Then there are only three vertices outside $N[x].$ We will use the fact that every vertex in $N_i(x)$
has at least one neighbor in $N_{i-1}(x)$ for $1\le i\le {\rm diam}(G).$
The proof consists of  a series of claims.

{\bf Claim 1.} ${\rm diam}(G)=3.$

Clearly ${\rm diam}(G)=ecc(x)\le 4.$ If $ecc(x)=4,$ let $x,r,s,p,q$ be a diametral path. Then $ecc(r)\le 3$
and $ecc(s)\le 3,$ implying that both $r$ and $s$ are central vertices, a contradiction. Thus ${\rm diam}(G)\le 3.$ On the other hand, it is impossible
that ${\rm diam}(G)=2,$ since otherwise ${\rm rad}(G)=1,$ implying that $\Delta(G)=n-1,$ a contradiction. Hence ${\rm diam}(G)=3.$

{\bf Claim 2.} The vertex $x$ has only one eccentric vertex, which is not a leaf.

If $|N_3(x)|=2,$ then $|N_2(x)|=1.$  Now the vertex in $N_2(x)$ and its neighbors in $N(x)$ are central vertices, a contradiction. Thus
$x$ has only one eccentric vertex, which we denote by $w.$ Let $N_2(x)=\{u,\,v\}.$ If $w$ is a leaf, without loss of generality, suppose $u$ is the neighbor
of $w.$ Since $ecc(w)\le 3,$ we deduce that $ecc(u)\le 2;$ i.e., $u$ is a central vertex. If $u$ and $v$ are adjacent, then every neighbor of $u$ in $N(x)$
is also a central vertex, a contradiction; if $u$ and $v$ are nonadjacent, then $d(u,\,v)=2,$ implying that $u$ and $v$ have a common neighbor $y$ in $N(x).$
But then $y$ is also a central vertex, a contradiction again.

Claim 2 shows that $N(w)=\{u,\,v\}.$

{\bf Claim 3.} $u$ and $v$ have at most one common neighbor in $N(x).$

This holds since every common neighbor of $u$ and $v$ in $N(x)$ is a central vertex.

{\bf Claim 4.} $u$ and $v$ are nonadjacent.

To the contrary, assume that $u$ and $v$ are adjacent. Then any neighbor of either $u$ or $v$ in $N(x)$ is a central vertex. It follows that
$u$ and $v$ have a common neighbor $y$ in $N(x)$ and $y$ is their only neighbor in $N(x).$ Now any vertex $z\in N(x)\setminus\{y\}$ must be adjacent to $y,$
since $d(z,\,w)\le 3.$ Consequently ${\rm deg}(y)=n-2>n-4=\Delta(G),$ a contradiction.

{\bf Claim 5.} Neither $u$ nor $v$ is a top vertex.

To the contrary, assume ${\rm deg}(u)=n-4.$ By Claim 4, $|N(u)\cap N(x)|=n-5.$ Let $z\in N(x)$ be the nonneighbor of $u.$  If $d(u,z)=2,$ then $u$ is a central vertex and $u,\,v$ have no common neighbor in $N(x).$ Hence $v$ is adjacent to $z.$ Let $y$ be a common neighbor of $u$ and $z.$ Then $y$ is also a central vertex, a contradiction. Hence $d(u,z)=3.$ The condition $d(z,w)\le 3$ implies that $z$ is adjacent to $v.$ If $u$ and $v$ have no common neighbor in $N(x),$ then  $G$ is self-centered, a contradiction. If $u$ and $v$ have a common neighbor in $N(x),$ then $x$ and $u$ are the only two vertices with degree $n-4,$
contradicting our assumption that $G$ has at least three top vertices.

Similarly we can prove that ${\rm deg}(v)<n-4.$

 {\bf Claim 6.} Each of $u$ and $v$ has at least two neighbors in $N(x).$

To the contrary, assume that $N(u)\cap N(x)=\{y\}.$ Then for any vertex $z\in N(x)\setminus\{y\},$ $z$ cannot be adjacent to both $y$ and $v$, since
otherwise $y$ and $z$ are central vertices. Considering $d(z, u)$ and $d(z, w)$ we deduce that $z$ is a peripheral vertex. If $y$ and $v$ are
 nonadjacent, then $G$ is self-centered, a contradiction. If $y$ and $v$ are adjacent, then $y$ is the central vertex. Since $d(z, w)\le 3,$ we obtain
$d(z, v)\le 2.$ It follows that $v$ is also a central vertex, a contradiction.

Similarly we can prove that $v$ has at least two neighbors in $N(x).$

{\bf Claim 7.} $G$ has at most $n-5$ top vertices and the extremal graph is unique.

By Claim 3 and Claim 6, $u$ has a neighbor $y$ in $N(x)$ that is nonadjacent to $v,$ and $v$ has a neighbor $z$ in $N(x)$ that is nonadjacent to $u.$
Note that if $f$ is a neighbor of $u$ in $N(x)$ and $g$ is a neighbor of $v$ in $N(x)$ with $f\neq g,$ then $f$ and $g$ are nonadjacent, since otherwise
$f$ and $g$ are central vertices. Using Claim 6 again we deduce that neither $y$ nor $z$ has maximum degree. Thus $G$ has at least the five vertices
$y,z,u,v,w$ with degrees less than $n-4.$ It follows that $G$ has at most $n-5$ top vertices.

Conversely, suppose $G$ has $n-5$ top vertices. Then the above analysis shows that (1) each of $u$ and $v$ has exactly two neighbors in $N(x);$
(2) $u$ and $v$ have exactly one common neighbor $h$ in $N(x);$ and (3) the closed neighborhood of every vertex in $N(x)\setminus \{y,z,h\}$ is equal to $N[x].$
Consequently $G$ is the graph obtained from the graph of order $7$ in Figure 4 by blowing up a non-central vertex of degree $3$ into $K_{n-6}.$ This completes
the proof. \hfill $\Box$

The extremal graph of order $10$ in Theorem 11 is depicted in Figure 6.
\vskip 3mm
\par
 \centerline{\includegraphics[width=3.2in]{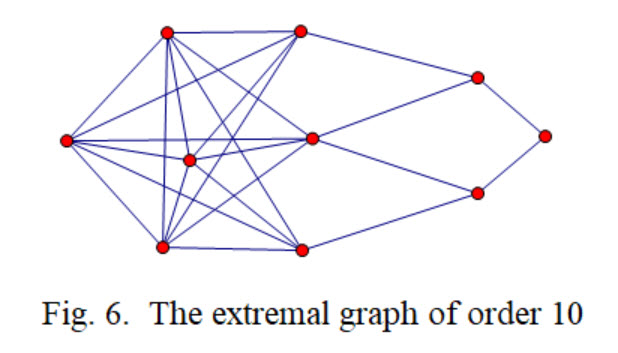}}
\par

\vskip 5mm
{\bf Acknowledgement.} This research  was supported by the NSFC grants 11671148 and 11771148 and Science and Technology Commission of Shanghai Municipality (STCSM) grant 18dz2271000.


\begin{thebibliography}{99}
\bibitem{1} B. Bollob\'{a}s and E. Szemer\'{e}di, Girth of sparse graphs, J. Graph Theory, 39(2002), no.3, 194-200.
\bibitem{2} F. Buckley and F. Harary, Distance in Graphs, Addison-Wesley Publishing Company, 1990.
\bibitem{3} F. Harary and R.Z. Norman, The dissimilarity characteristic of Husimi trees, Ann. of Math., 58(1953), no.1, 134-141.
\bibitem{4} S. Klav\v{z}ar, K.P. Narayankar and H.B. Walikar, Almost self-centered graphs, Acta Math. Sin. (Engl. Ser.), 27(2011), no.12, 2343-2350.
\bibitem{5} S. Klav\v{z}ar, K.P. Narayankar, H.B. Walikar and S.B. Lokesh, Almost peripheral graphs, Taiwanese J. Math., 18(2014), no.2, 463-471.
\bibitem{6} L. Lesniak, Eccentric sequences in graphs, Period. Math. Hungar., 6(1975), no.4, 287-293.
\bibitem{7} D.B. West, Introduction to Graph Theory, Prentice Hall, Inc., 1996.
\end{thebibliography}
\end{document}